\documentclass[12pt]{article}
\usepackage{amsmath}
\usepackage{amssymb}
\usepackage{amsmath,amssymb,amsbsy,amsfonts,amsthm,latexsym,
amsopn,amstext,amsxtra,euscript,amscd}

\begin{document}

\def\A{\mathbb{A}}
\def\B{\mathbf{B}}
\def \C{\mathbb{C}}
\def \F{\mathbb{F}}
\def \K{\mathbb{K}}

\def \Z{\mathbb{Z}}
\def \P{\mathbb{P}}
\def \R{\mathbb{R}}
\def \Q{\mathbb{Q}}
\def \N{\mathbb{N}}
\def \Z{\mathbb{Z}}

\def\B{\mathcal B}
\def\e{\varepsilon}

\def\cA{{\mathcal A}}
\def\cB{{\mathcal B}}
\def\cC{{\mathcal C}}
\def\cD{{\mathcal D}}
\def\cE{{\mathcal E}}
\def\cF{{\mathcal F}}
\def\cG{{\mathcal G}}
\def\cH{{\mathcal H}}
\def\cI{{\mathcal I}}
\def\cJ{{\mathcal J}}
\def\cK{{\mathcal K}}
\def\cL{{\mathcal L}}
\def\cM{{\mathcal M}}
\def\cN{{\mathcal N}}
\def\cO{{\mathcal O}}
\def\cP{{\mathcal P}}
\def\cQ{{\mathcal Q}}
\def\cR{{\mathcal R}}
\def\cS{{\mathcal S}}
\def\cT{{\mathcal T}}
\def\cU{{\mathcal U}}
\def\cV{{\mathcal V}}
\def\cW{{\mathcal W}}
\def\cX{{\mathcal X}}
\def\cY{{\mathcal Y}}
\def\cZ{{\mathcal Z}}

\def\f{\frac{|\A||B|}{|G|}}
\def\AB{|\A\cap B|}
\def \Fq{\F_q}
\def \Fqn{\F_{q^n}}

\def\({\left(}
\def\){\right)}
\def\fl#1{\left\lfloor#1\right\rfloor}
\def\rf#1{\left\lceil#1\right\rceil}
\def\Res{{\mathrm{Res}}}

\newcommand{\comm}[1]{\marginpar{
\vskip-\baselineskip \raggedright\footnotesize
\itshape\hrule\smallskip#1\par\smallskip\hrule}}

\newtheorem{lem}{Lemma}
\newtheorem{lemma}[lem]{Lemma}
\newtheorem{prop}{Proposition}
\newtheorem{proposition}[prop]{Proposition }
\newtheorem{thm}{Theorem}
\newtheorem{theorem}[thm]{Theorem}
\newtheorem{cor}{Corollary}
\newtheorem{corollary}[cor]{Corollary}
\newtheorem{prob}{Problem}
\newtheorem{problem}[prob]{Problem}
\newtheorem{ques}{Question}
\newtheorem{question}[ques]{Question}
\newtheorem{rem}{Remark}

\title{A note on $n!$ modulo $p$}

\author{
{\sc M. Z. Garaev} and {\sc J. Hern\'andez}}

\date{}

\maketitle

\begin{abstract}
Let $p$ be a prime, $\varepsilon>0$ and $0<L+1<L+N < p$. We prove
that if $p^{1/2+\varepsilon}< N <p^{1-\varepsilon}$, then
$$
\#\{n!\!\!\! \pmod p;\,\, L+1\le n\le L+N\} > c (N\log N)^{1/2},\,\, c=c(\varepsilon)>0.
$$
We use this bound to show that any $\lambda\not\equiv 0\pmod p$ can
be represented in the form $\lambda \equiv n_1!\ldots n_7!\pmod p$,
where $n_i=o(p^{11/12})$. This slightly refines the previously known
range for $n_i$.
\end{abstract}

\section{Introduction}

In what follows, $p$ is a large prime number. For integers $L$ and
$N$ with
$$
0<L+1<L+N < p
$$
we consider the set
$$
\cA(L, N)=\Bigl\{n!\!\! \pmod p;\,\, L+1\le n\le L+N\Bigr\}.
$$
From the observation
\begin{equation}
\label{eqn:int inclusion}
\{1\}\cup\{L+2,\ldots, L+N\}\pmod p \subset \frac{\cA(L, N)}{\cA(L, N)}
\end{equation}
it follows that
$$
|\cA(L, N)|\ge N^{1/2}.
$$
In particular, we trivially have $|\cA(0,p-1)|\ge (p-1)^{1/2}$. From
the result of Garc\'ia~\cite{VC2} on the cardinality of product of
two factorials modulo $p$ it follows that  $|\cA(0,p-1)|> cp^{1/2}$
for any constant $c<\sqrt{\frac{41}{24}}$ and sufficiently large
prime $p$.  The conjecture is that
$|\cA(0,p)|$  asymptotically behaves like $(1-e^{-1})p$,
see~\cite{CVZ} and~\cite{Guy}.

Improving on the trivial bound, Klurman and Munsch~\cite{KM} proved the bound
\begin{equation}
\label{eqn:constanta c}
|\cA(L, N) |\ge c N^{1/2}
\end{equation}
with $c=\sqrt{\frac{3}{2}}$ and $p^{1/4+\varepsilon}<N<p$. We remark that the condition
$N>p^{1/4+\varepsilon}$ can be relaxed if one combines~\eqref{eqn:int inclusion}
with~\cite[Theorem 2]{CillGar} (see, also, the work~\cite{CCGHSZ} on congruences
with variables from short intervals).

In the present note, using a consequence of Bombieri's bound on
exponential sums over algebraic curves, we show that if
$p^{1/2+\varepsilon} < N = o(p)$, then the constant $c$
in~\eqref{eqn:constanta c} can be taken arbitrarily large. We then
apply this result to the problem of representability of residue
classes as a product of seven factorials with small variables.
\begin{theorem}
\label{thm:n!} Let $p^{1/2+\varepsilon}< N <0.1p$. Then
$$
\Bigl|\frac{\cA(L, N)}{\cA(L, N)}\Bigr|>c_0 N\log(p/N)
$$
for some $c_0=c_0(\varepsilon)>0.$
\end{theorem}

From Theorem~\ref{thm:n!} it follows, in particular, that for
$p^{1/2+\varepsilon}< N <0.1p$ we have the bound
$$
|\cA(L, N)|>c_0 (N\log(p/N))^{1/2}
$$
for some $c_0=c_0(\varepsilon)>0.$

Garaev, Luca and Shparlinski~\cite{GLSh} proved that any
$\lambda\not\equiv 0\pmod p$ can be represented in the form
$$
\prod_{i=1}^7 n_i!\equiv \lambda\pmod p,
$$
where $n_i\ll p^{11/12+\varepsilon}$.
Garcia~\cite{VC1} improved this condition to $n_i\ll p^{11/12}$.
Using Theorem~\ref{thm:n!} we can slightly improve this as follows.

\begin{theorem}
\label{thm:n!7} Any $\lambda\not\equiv 0\pmod p$ can be represented
in the form
$$
\prod_{i=1}^7 n_i!\equiv \lambda\pmod p,
$$
where the positive integers $n_1,...,n_7$ satisfy
$$
\max\{n_i | i=1,...,7\} \ll\frac{p^{11/12}}{(\log p)^{1/2}}.
$$
\end{theorem}

\bigskip

\section{Lemmas}

We need the following special case of the results of Bombieri~\cite[Theorem 6]{Bom} and Chalk and Smith~\cite[Theorem 2]{ChS}.
As usual, $\F_p$ denotes the field of residue classes modulo $p$.

\begin{lemma}
\label{lem:Bom} Let $(b_1,b_2)\in \F_p\times \F_p$ be nonzero and
$f(x,y)\in \F_p[x,y]$ be a polynomial of degree $d\ge 1$ with the
following property: there is no $c\in \F_p$ for which the polynomial
$f(x,y)$ is divisible by $b_1x+b_2y+c$. Then
$$
\Bigl|\sum_{f(x,y)=0}e^{2\pi i(b_1x+b_2y)/p}\Bigr|\le 2d^2 p^{1/2}.
$$
\end{lemma}

We remark that the factor $2$ on the right hand side can be removed, but it is not essential in our application.

The following lemma is due to Ruzsa. It will be used in the proof of Theorem~\ref{thm:n!7}.

\begin{lemma}
\label{lem:RUZ} For any finite subsets $X,Y,Z$ of an abelian group
we have
$$
|X-Y|\le \frac{|X+Z||Z+Y|}{|Z|}.
$$
\end{lemma}

In the proof of Theorem~\ref{thm:n!7} we will also need the following estimate of character sums with
factorials from the work of Garc\'ia~\cite{VC1}.

\begin{lemma}
\label{lem:VG} For any positive integer $N$ the following bound
holds:
$$
\max_{\chi\not=\chi_0}\Bigl|\sum_{n\le N}\sum_{m\le N}\chi((n+m)!)\Bigr|\ll N^{7/4} p^{1/8}.
$$
\end{lemma}

\section{Proof of Theorem~\ref{thm:n!}}

We can assume that $p/N$ is sufficiently large in terms of
$\varepsilon$. Let
$$
M=\lfloor\min\{p^{0.1\varepsilon},
(p/N)^{0.1}\}\rfloor
$$
For a positive integer $j\le M$ we define the set
$$
X_j=\Bigl\{\prod_{i=1}^j(x+L+i)\pmod p;\quad 1\le x < 0.6N\Bigr\}.
$$
Since the polynomial $\prod_{i=1}^j(x+L+i)$ has degree $j$, we have
that
\begin{equation}
\label{eqn: X_j> N/j}
|X_j|\ge \frac{N}{2j}.
\end{equation}
Let us prove that for any $j\ge 2$ the following bound holds:
$$
\#\{X_j\setminus(X_1\cup\ldots\cup X_{j-1})\}\ge \frac{N}{3j}.
$$

Note that
\begin{eqnarray*}
\begin{split}
&\#\{X_j\setminus(X_1 \cup \ldots\cup X_{j-1})\}  \\ &\qquad \qquad = \#\{X_j\setminus((X_j\cap X_1)\cup\ldots\cup (X_j\cap X_{j-1}))\}\\
 &\qquad \qquad \qquad \qquad \ge |X_j|-|X_j\cap X_1|-\ldots |X_j\cap X_{j-1}|.
\end{split}
\end{eqnarray*}
Therefore, in view of~\eqref{eqn: X_j> N/j} we get
\begin{equation}
\label{eqn:X_j setminus}
\#\{X_j\setminus(X_1\cup\ldots\cup X_{j-1})\}\ge  \frac{N}{2j} - |X_j\cap X_1|-\ldots -|X_j\cap X_{j-1}|.
\end{equation}
We shall obtain upper bound for the cardinality $|X_j\cap X_k|$ for
$1\le k\le j-1$. Let $J(j,k)$ be the number of solutions of the
congruence
$$
\prod_{i=1}^j(x+L+j)\equiv \prod_{i=1}^k(y+L+i)\pmod p,\quad 1\le x,y < 0.6N.
$$
Clearly,
\begin{equation}
\label{eqn:X_j J(j,k)}
|X_j\cap X_k|\le J(j,k).
\end{equation}
Denote
$$
f(x,y)= \prod_{i=1}^j(x+L+i) - \prod_{i=1}^k(y+L+i)\in \F_p[x,y].
$$
Following standard arguments,  we write $J(j,k)$ in the form
\begin{eqnarray*}
J(j,k)&=&\sum_{\substack{x\le 0.6N, \,y\le 0.6N\\ f(x,y)=0}}1 \\
&\ge& \frac{1}{p^2}\sum_{b_1=0}^{p-1}\sum_{b_2=0}^{p-1}\sum_{u<0.6N}\sum_{v<0.6N}\sum_{f(x,y)=0}e^{2\pi i (b_1(x-u)+b_2(y-v))/p}.
\end{eqnarray*}
From the trivial bound we have that the number of solutions of the
equation
$$
f(x,y)=0,\quad (x,y)\in \F_p\times \F_p
$$
is not greater, than $jp$. We also recall the elementary estimates
$$
\sum_{b=0}^{p-1}\Bigl|\sum_{z<0.6N} e^{2\pi i buz/p}\Bigr|<p\log p.
$$
 Thus, separating the term
that corresponds to $b_1=b_2=0$, we obtain
$$
J(j,k)\le \frac{jN^2}{p} +
(\log p)^2\max_{(b_1,b_2)}\Bigl|\sum_{f(x,y)=0}e^{2\pi i (b_1x+b_2y)/p}\Bigr|,
$$
where the maximum is taken over the integers $0\le b_1,b_2\le p-1$
such that $(b_1,b_2)\not=(0,0)$. Since $j>k\ge 1$, for any $a_1,a_2,a_3\in \F_p$ the polynomials
$f(X,a_1X+a_2)$ and $f(a_3,X)$ are polynomials of degrees $j$ and $k$ in $\F_p[X]$. Therefore,
$f(x,y)$ is not divisible by $b_1x+b_2y+c$ in $\F_p[x,y]$ and thus satisfies the condition of Lemma~\ref{lem:Bom}. Hence,
applying Lemma~\ref{lem:Bom} and taking into account that $j\le M$, we
get
$$
J(j,k)\le \frac{jN^2}{p} +
O((\log p)^2 j^2 p^{1/2})\le \frac{N}{6j^2}.
$$
This bound and~\eqref{eqn:X_j J(j,k)} together with~\eqref{eqn:X_j
setminus} implies that
$$
\#\{X_j\setminus(X_1\cup\ldots\cup X_{j-1})\}\ge \frac{N}{2j}-\frac{(j-1)N}{6j^2}\ge \frac{N}{3j}.
$$

Now we  observe that
$$
X_j\!\pmod p\subset \frac{\cA(L, N)}{\cA(L, N)},\quad j=1,2,\ldots, M.
$$
Hence
\begin{eqnarray*}
\Bigl|\frac{\cA(L, N)}{\cA(L, N)}\Bigr| &\ge& \#\{X_1\cup X_2\cup \ldots X_m\}\\
&=& |X_1|+\sum_{j=2}^m\#\{X_j\setminus(X_1\cup\ldots\cup X_{j-1})\} \\
&\ge& \sum_{j=1}^M\frac{N}{3j}\gg N\log M\gg N\log(p/N)
\end{eqnarray*}
and the result follows.

\section{Proof of Theorem 2}

Let $p^{0.51}< N < p^{0.99}$. For the brevity, denote $\cA=\cA(0,
N)$.  By Theorem~\ref{thm:n!} we have
$$
\Bigl|\frac{\cA}{\cA}\Bigr|\gg N\log p,\quad |\cA|\gg (N\log p)^{1/2}.
$$
Application of Lemma~\ref{lem:RUZ} in the multiplicative form gives
the bound
$$
\Bigl|\frac{\cA}{\cA}\Bigr| \le \frac{|\cA \cA|^2}{|\cA|}
$$
Hence,
\begin{equation}
\label{eqn:AA>}
|\cA\cA|\ge c_1(N \log p)^{3/4}
\end{equation}
for some absolute constant $c_1>0$.

Denote $I=\{1,2,\ldots, N\}$. Let $J$ be the number of solutions of
the congruence
$$
(n_1+m_1)!(n_2+m_2)!(n_3+m_3)! x y \equiv \lambda \pmod p,
$$
in variables $n_1,n_2,m_1,m_2,x,y$ satisfying
$$
n_1,n_2,m_1,m_2\in I,\quad x,y\in \cA\cA.
$$
To prove Theorem~\ref{thm:n!7} it suffices to show that there is a
constant $C>0$ such that $J>0$ for $N =\lfloor Cp^{11/12}\rfloor$. We
express $J$ via character sums and get
$$
J=\frac{1}{p-1}\sum_{\chi}\sum_{n_1,n_2,m_1,m_2\in I}\,\,\sum_{x,y\in \cA\cA}\chi((n_1+m_1)!(n_2+m_2)!(n_3+m_3)! x y)\chi(\lambda^{-1}).
$$
Separating the term that corresponds to the principal character
$\chi=\chi_0$ and following the standard argument we obtain
$$
J\ge \frac{N^6 |\cA\cA|^2}{p-1}-\frac{1}{p-1}\sum_{\chi\not=\chi_0}\Bigl|\sum_{n,m\in I}\chi((n+m)!)\Bigr|^3\Bigl|\sum_{x\in \cA\cA}\chi(x)\Bigr|^2.
$$
Application Lemma~\ref{lem:VG} and the identity
$$
\frac{1}{p-1}\sum_{\chi}\Bigl|\sum_{x\in \cA\cA}\chi(x)\Bigr|^2=|\cA\cA|,
$$
gives
$$
J\ge \frac{N^6 |\cA\cA|^2}{p-1}-c_2 N^{21/4}p^{3/8}|\cA\cA|,
$$
where $c_2>0$ is an absolute constant. Using~\eqref{eqn:AA>} we
obtain
\begin{eqnarray*}
J&\ge&  \frac{N^{21/4}|\cA \cA|}{p-1}\Bigl(|\cA\cA| N^{3/4}-c_2 p^{11/8}\Bigr)\\
&\ge& \frac{N^{21/4}|\cA\cA|}{p-1}\Bigl(c_1N^{3/2}(\log p)^{3/4}-c_2 p^{11/8}\Bigr).
\end{eqnarray*}
Hence, taking $N=\lceil 2(c_2/c_1)^{2/3}p^{11/18}/(\log p)^{1/2}\rceil$,
we get $J>0$, which finishes the proof of our theorem.

\section{Remarks}

In the proof of Theorem~\ref{thm:n!7} we used the fact that for $N<
p^{1-\varepsilon}$ one has the bound
$$
|\cA(0,N) \cA(0,N)|\gg (N\log N)^{3/4}.
$$
We note that this bound can significantly be improved for small
values of $N$. For example, let $N<p^{1/2}$. For any positive
integers $n,m\le N$ we have
$$
\frac{n}{m} \pmod p\subset \frac{\cA(0,N) \cA(0,N)}{\cA(0,N) \cA(0,N)}.
$$
Note that in the range $n,m<p^{1/2}$ for distinct rational numbers
$n/m$ correspond distinct residue classes $n/m\pmod p$. Therefore,
\begin{eqnarray*}
\Bigl|\frac{\cA(0,N) \cA(0,N)}{\cA(0,N) \cA(0,N)}\Bigr|&\ge& \#\Bigl\{\frac{n}{m};\,\, n,m\in [1,N]\cap\Z,\, \gcd(n,m)=1 \Bigr\}\\
&=&\Bigl(\frac{6}{\pi^2}+o(1)\Bigr)N^2
\end{eqnarray*}
as $N\to\infty$. Thus, in the range $N<p^{1/2}$ we have $|\cA(0,
N)\cA(0,N)|\gg N.$

\bigskip

{\bf Acknowledgement}. M. Z. Garaev was supported by the sabbatical grant from PASPA-DGAPA-UNAM.

Address of the authors:\\

M.~Z.~Garaev, Centro de Ciencias Matem\'{a}ticas, Universidad
Nacional Aut\'onoma de M\'{e}xico, C.P. 58089, Morelia,
Michoac\'{a}n, M\'{e}xico,

Email: {\tt garaev@matmor.unam.mx}

\vspace{1cm}

J. Hern\'andez, Centro de Ciencias Matem\'{a}ticas, Universidad
Nacional Aut\'onoma de M\'{e}xico, C.P. 58089, Morelia,
Michoac\'{a}n, M\'{e}xico,

Email: {\tt stgo@matmor.unam.mx}

\end{document}